\documentclass[11pt,leqno]{amsart}
\usepackage{amssymb,verbatim,enumerate,ifthen}
\usepackage[mathscr]{eucal}
\usepackage[utf8]{inputenc}
\usepackage[T1]{fontenc}
\usepackage{cite}
\def\N{\mathbb{N}}
\def\R{\mathbb{R}}

\def\Z{\mathbb{Z}}

\def\G{\mathscr{G}}

\def\L{\mathscr{L}}

\def\A{\mathscr{A}}

\newtheorem{theorem}{Theorem}[section]
\newtheorem*{theorem*}{Theorem}
\def\Thm#1#2{\ifthenelse{\equal{#1}{*}}{\begin{theorem*}#2\end{theorem*}}
             {\begin{theorem}\label{T#1}#2\end{theorem}}}
\newtheorem{Atheorem}{Theorem}

\def\thm#1{Theorem~\ref{T#1}}
\newtheorem{proposition}[theorem]{Proposition}
\newtheorem*{proposition*}{Proposition}
\def\Prp#1#2{\ifthenelse{\equal{#1}{*}}{\begin{proposition*}#2\end{proposition*}}
{\begin{proposition}\label{P#1}#2\end{proposition}}}
\def\prp#1{Proposition~\ref{P#1}}

\newtheorem{corollary}[theorem]{Corollary}
\newtheorem*{corollary*}{Corollary}
\def\Cor#1#2{\ifthenelse{\equal{#1}{*}}{\begin{corollary*}#2\end{corollary*}}
             {\begin{corollary}\label{C#1}#2\end{corollary}}}

\newtheorem{lemma}[theorem]{Lemma}
\newtheorem*{lemma*}{Lemma}
\def\Lem#1#2{\ifthenelse{\equal{#1}{*}}{\begin{lemma*}#2\end{lemma*}}
             {\begin{lemma}\label{L#1}#2\end{lemma}}}

\theoremstyle{definition}
\newtheorem{remark}[theorem]{Remark}
\newtheorem*{remark*}{Remark}
\def\Rem#1#2{\ifthenelse{\equal{#1}{*}}{\begin{remark}\rm #2\end{remark}}
             {\begin{remark}\label{R#1}\rm #2\end{remark}}}

\newtheorem{example}[theorem]{Example}
\newtheorem*{example*}{Example}
\def\Exa#1#2{\ifthenelse{\equal{#1}{*}}{\begin{example*}\rm #2\end{example*}}
             {\begin{example}\label{Ex#1}\rm #2\end{example}}}

\def\eq#1{{\rm(\ref{E#1})}}
\def\Eq#1#2{\ifthenelse{\equal{#1}{*}}
  {\begin{equation*}\begin{aligned}#2\end{aligned}\end{equation*}}
  {\begin{equation}\begin{aligned}\label{E#1}#2\end{aligned}\end{equation}}}

\begin{document}
\begin{flushright}
\end{flushright}
\vspace{5mm}

\date{\today}

\title[On characterizations, Decompositions, and Stability of Convex Sequences]
{On characterizations, Decompositions, and Stability of Convex Sequences}

\author[A. R. Goswami]{Angshuman R. Goswami}
\address[A. R. Goswami]{Department of Mathematics, University of Pannonia,
H-8200 Veszprém, Hungary}
\email{goswami.angshuman.robin@mik.uni-pannon.hu}



\begin{abstract}
This paper introduces new characterizations, decomposition theorems, and stability results for convex sequences. We show that a sequence is convex precisely when its epigraph satisfies a midpoint convexity condition, thereby connecting discrete and geometric notions of convexity. A decomposition result proves that any sequence can be written as the difference of two convex sequences, with generalizations to higher-order convexity. We construct nontrivial convex minorants for bounded-below sequences and establish a Hyers-Ulam-type stability theorem showing that any approximately convex sequence can be uniformly approximated by a genuine convex sequence without significantly altering its values. Finally, for a concave 
sequence, we characterize those subsequences that are convex in it by providing slope inequalities and monotone auxiliary sequences. We explore the interplay among convexity, subadditivity, and periodically indexed subsequences.  
\end{abstract}
\maketitle
\section*{Introduction}
Throughout this paper $\N$, $\Z$, $\R$, and $\R_+$ are used to denote the sets of natural, integer, real, and non-negative real numbers, respectively. For a sequence $\big(u_n\big)_{n=0}^{\infty}$, the symbol $\{u_n\}$ represents its range set. For greater generality, our analysis is primarily carried out for infinite sequences. Nevertheless, the results admit natural finite analogues, which can be derived without substantial changes to the arguments.\\

A sequence $\big(u_n\big)_{n=0}^{\infty}$ is \textit{convex} if for all $i,j\in\N$ with $i\leq j$, the following functional inequality is satisfied
\Eq{9910}{
u_i-u_{i-1}\leq u_{j+1}-u_{j}.
} 
In other words, a sequence $\big(u_n\big)_{n=0}^{\infty}$ possesses convexity if the following discrete functional inequality holds
\Eq{01}{
2u_n\leq u_{n-1}+u_{n+1}\qquad \mbox{for all}\qquad n\in\N.
}
If the reverse inequality holds, we call it a \textit{concave} sequence. Clearly, the class of convex sequences contains some of the well-known classes of sequences such as arithmetic, geometric, partitioning, factorial, Fibonacci, etc. \\

The first mention of convex sequences appears in the book \textit{Analytic Inequalities} (see \cite{Mitrinovic}) in 1970. As the discrete version of convex functions, convex sequences also possess many interesting properties. Over the last few decades, researchers have investigated various structural characteristics, characterizations, generalizations, and related inequalities of convex sequences. Some investigations towards higher-order and approximate convex sequences have also been carried out. Moreover, several algorithms based on convex sequences have been developed, demonstrating their potential applications beyond mathematics. Most of these contributions can be found in \cite{Murota, Murotaa, Shioura, Shiouraa, Goswami2025a, Goswami2025b, Goswami2025c, pecaric, Essen, Jimenez, Mercer,Krasniqi, Sofonea}.\\

In this paper, we contribute to the rich history of this topic and fill some gaps by addressing several natural but challenging questions. The paper is organized into three sections.
\\

In the first section, we provide a characterization of convex sequences. This investigation is motivated by the close connections of epigraphs and subadditivity with discrete convexity.  Our initial structural analysis also addresses the decomposition of arbitrary sequences, a notion that draws on the classical theorem on the decomposition of functions of bounded variation. This result, due to Camille Jordan (1881), arose from his study of the convergence of Fourier series and is formulated as follows (see \cite{Jordan})
\Thm{*}{\normalfont{\textbf{[Jordan (see \cite{Jordan})]}} In a compact interval, any function of bounded variation can be represented as the difference of two monotone(increasing) functions.
}
Analogous to this theorem, the following result was demonstrated in \cite{Goswami2025c}.
\Thm{*}{\normalfont{\textbf{[(see \cite{Goswami2025c})]}} Any sequence $\big(u_n\big)_{n=0}^{\infty}$ can be re-written as the difference of two monotonically increasing sequences $\big(v_n\big)_{n=0}^{\infty}$ and $\big(w_n\big)_{n=0}^{\infty}$; where the sequence $\big(v_n\big)_{n=0}^{\infty}$ act as a monotone majorant of $\big(u_n\big)_{n=0}^{\infty}$ and the sequence $\big(w_n\big)_{n=0}^{\infty}$ possesses non-negativity besides monotonicity.}

In this paper, we present a way to express any sequence as the difference between two convex sequences and discuss the possibility of generalizing this result to higher-order sequential convexity.\\ 

In the second section, we discuss the possibilities to formulate a non-trivial convex minorant for any given sequence. We also present a sandwich-type result by discussing the constraints two distinct non-overlapping sequences need to satisfy for a convex sequence to 
exist between them. Using these two results, we prove that if a sequence preserves the convexity property in an approximate sense, then it is possible to bring back the ideal convex structure of it without significantly altering any of the sequential values. In classical function theory, such results are referred to as Hyers-Ulam-type stability results, named after the pioneering and influential work of the mathematicians D. H. Hyers and S. M. Ulam. In view of a single variable function, some of the main results obtained by them can be stated as follows
\Thm{*}{\normalfont{\textbf{[Hyers (see \cite{Ulam})]}} Let $\varepsilon>0$ be fixed. If $f:\R\to\R$ satisfies the following  functional inequality
\Eq{*}{
|f(x+y)-f(x)-f(y)|\leq \varepsilon \quad \mbox{for all}\quad x,y\in \R,
}
then corresponding to $f$, there exists an additive function $\A:\R\to\R$ such that the inequality holds
\Eq{*}{
|f(x)-A(x)|\leq \varepsilon \quad \mbox{for all}\quad x\in\R.
}
}
\Thm{*}{\normalfont{\textbf{[Hyers, Ulam (see \cite{Hyers})]}} Let $\varepsilon>0$ be fixed. If $f:\R\to\R$ satisfies the following functional inequality
\Eq{*}{
f(tx+(1-t)y)\leq tf(x)+(1-t)f(y)+\varepsilon \quad \mbox{for all}\quad x,y\in R\quad \mbox{and} \quad t\in[0,1],
}
then corresponding to $f$, there exists a convex function $g:\R\to\R$ such that the inequality holds
\Eq{*}{
|f(x)-g(x)|\leq \dfrac{\varepsilon}{2} \quad \mbox{for all}\quad x\in\R.
}
}
These earlier investigations, related results, and their implementation in discrete settings in later stage can be found in the articles \cite{Hyers, Becken, Ulam, Pales, Goswami2025a,Goswami2025b,Goswami2025c} which basically demonstrate that if a function $f$ satisfies a given equation or inequality with a sufficiently small error $\varepsilon$, then there may exists an exact solution $g$ whose distance from $f$ is controlled by $\varepsilon$.\\

As a discrete version of the convexity result, we prove that for a fixed $\varepsilon>0$, if a sequence $\big(u_n\big)_{n=0}^{\infty}$ is bounded from below, and satisfies the following discrete functional inequality 
\Eq{*}{
u_n\leq\sum_{i=1}^{n-1}p_i\,u_i+q\,u_n+\sum_{i=n+1}^{\infty}r_i\,u_i+\varepsilon\qquad\mbox{with}\qquad \sum_{i=1}^{n-1} p_i+q+\sum_{i=n+1}^\infty r_i=1\,\,;
}
where $p_i$'s, q and $r_i$'s are proper non-negative dyadic rational numbers such that either some $p_i$'s and some $r_i$'s simultaneously remain positive or else $q=1$. Then there exists a convex sequence 
$\big(v_n\big)_{n=0}^{\infty}$ such that the norm inequality $\|u_n-v_n\|_{\infty}\leq \varepsilon/2$ holds. We also propose a converse assertion of this result.\\

In the final section of this paper, we perform some analysis that may give better insights into Erd\"os-Szekeres combinatorial results in some particular cases. Their monotone subsequence result, one of the highly celebrated results in combinatorial mathematics can be stated as follows
\Thm{*}{\normalfont{\textbf{[Erd\"os-Szekeres (see \cite{ErdosSzekeres})]}}
Let $n\in\N$. Every real sequence of length $n^2+1$ contains a monotone (either increasing or decreasing) subsequence of length at least$n+1$.
}

However, this combinatorial breakthrough has several limitations. In particular, it does not reveal any explicit procedure for constructing a monotone subsequence from an arbitrary given sequence. Moreover, even under the cardinality assumptions, the monotone subsequence theorem does not identify any specific properties of the original sequence that determine whether it must contain an increasing subsequence, a decreasing subsequence, or both. \\

It is also worth mentioning that while investigating a discrete geometry problem, Erd\"os and Szekeres obtained a result describing the length of the largest possible convex (or concave) subsequence, provided no three consecutive terms are equal. However, they did not use the usual terminology of sequential convexity at that time. Their combinatorial result on convexity is often referred to as the cup-cap theorem (see \cite{ErdosSzekeres1960}).\\

A monotonically increasing sequence contains no decreasing subsequence. Interestingly, this is not the case for convexity. There exist convex sequences that contain concave subsequences, and vice versa.
For example, the square root sequence 
$\big(\sqrt{n}\big)_{n=0}^{\infty}$ is concave. Now, if we choose indices $n_k's\in\N\cup\{0\}$ such that it satisfies the following inequality condition
\Eq{*}{ 
 n_k\leq\dfrac{\A\big(n_{_{k-1}},n_{_{k+1}}\big)+\G\big(n_{_{k-1}},n_{_{k+1}}\big)}{2}\qquad
\mbox{where}\qquad
\A\big(n_{_{k-1}},n_{_{k+1}}\big):&=\dfrac{n_{_{k-1}}+n_{_{k+1}}}{2}\\
\mbox{and}\qquad
 \G\big(n_{_{k-1}},n_{_{k+1}}\big):&=\sqrt{n_{_{k-1}}n_{_{k+1}}},
}
then the resultant subsequence of $\Big(\sqrt{n_{k}}\Big)_{k=0}^{\infty}$ turns to be convex. However, there exist concave sequences that contain no convex subsequence, and vice versa. One such example is the sequence $\big(2^n\big)_{n=0}^{\infty}$. It possesses convexity, and it is not possible to extract any concave sequence from it. We provide a characterization that describes the necessary and sufficient conditions under which a concave sequence carries a convex subsequence. A similar investigation for subadditive sequences is also carried out. Moreover, we show that for any increasing convex sequence, the subsequence obtained by selecting terms at indices that themselves form a convex sequence in $\N$ is also convex.\\

Combining such findings with Erd\"os-Szekeres-type combinatorial results along with some conditions may provide more insights towards the study of some specific sequence-subsequence relationships. More precisely, such formulations give both combinatorial bounds and required characterizations under which a sequence $\big(u_i)$ with property $p$ and length $\varrho(n)$ contains a subsequence $\big(u_{i_{_{k}}}\big)$ of length $n$ that possesses characteristics $q$.
\\

We start our investigation with some basic structural properties of convex sequences.

\section{On Characterization and Decomposition Results}
To present our first characterization, we have to recall the definition of the epigraph and an interpolation result for convex sequences. For a function $f:I(\subseteq\R)\to\R$, the \textit{epigraph} of $f$ or $epi(f)$ is defined as follows
\Eq{*}{
epi\big(f\big):=\Big\{(x,y)\,\,\big|\,\,  x\in I\,,\, f(x)\leq y\Big\}\subseteq I\times\R.
}
It is well known that a function $f:I(\subseteq\R)\to\R$ possesses convexity if and only if $epi(f)$ is a convex set. We want to establish a similar result involving a sequence and its epigraph.
Mathematically, the epigraph of the sequence $\big(u_n\big)_{n=0}^{\infty}$ can be  expressed as 
\Eq{*}{
epi\big(u_n\big):=\Big\{(n,y)\,\,\big|\,\,  n\in\N\cup\{0\}\,,\, u_n\leq y\Big\}\subseteq\big(\N\cup\{0\}\big)\times\R.
}
Also, we need the following theorem (see \cite{Goswami2025c}) which states that corresponding to every convex sequence there exists a convex function that perfectly interpolates the sequence.
\Thm{15}{Let $\big(u_n\big)_{n=0}^{\infty}$ be a convex sequence. Then there exists a continuous almost everywhere differentiable convex function $f:\R_{+}\to\R$ such that $f(n)=u_{n}$ for all $n\in\N\cup\{0\}.$\\
Conversely, if $f$ is convex in $\R_+$, then the sequence $\big(f(n)\big)_{n=0}^{\infty}$ is also convex.}
In the proof of the theorem, it was demonstrated that for each convex sequence $\big(u_n\big)_{n=0}^{\infty}$, the function $f:\R_+\to\R$ formed by joining line segments between each pair of consecutive points possesses convexity. 
\Prp{01}{
A sequence $\big(u_n\big)_{n=0}^{\infty}$ is convex if and only if for all $n_1,n_2\in\N\cup\{0\}$ with $n_1 \equiv n_2 \pmod{2}$, the inclusion $\bigg(\dfrac{n_1+n_2}{2}\,,\,\dfrac{u_{n_1}+u_{n_2}}{2}\bigg)\in epi \big(u_n\big)$ is satisfied.
}
\begin{proof}
To establish the first part of the assertion, we assume 
$\big(u_n\big)_{n=0}^{\infty}$ is convex. Thus \thm{15} holds, and results in the inclusion $epi\big(u_n\big)\subset epi(f).$ Now, we choose $n_1,n_2\in\N\cup\{0\}\subseteq\R_+$ arbitrarily such that the condition $n_1 \equiv n_2 \pmod{2}$ holds. From the definition of epigraphs, we can conclude the following
\Eq{*}{
\bigg(\dfrac{n_1+n_2}{2}\,,\,\dfrac{f({n_1})+f({n_2})}{2}\bigg)=\bigg(\dfrac{n_1+n_2}{2}\,,\,\dfrac{u_{n_1}+u_{n_2}}{2}\bigg)\in epi(f).
}
This implies $f(n)=u_n\leq \dfrac{u_{n_1}+u_{n_2}}{2}$ where $n=\dfrac{n_1+n_2}{2}\in\N$. It establishes the validity of the following inclusion under the mentioned conditions
\Eq{17}{
\bigg(\dfrac{n_1+n_2}{2}\,,\,\dfrac{u_{n_1}+u_{n_2}}{2}\bigg)\in epi \big(u_n\big)\quad \mbox{where}\quad n_1,n_2\in\N\cup\{0\}\quad\mbox{with}\quad n_1 \equiv n_2 \pmod{2}.
}
To prove the converse part, we assume that the inclusion \eq{17} holds. Now, replacing $n_1=n-1$ and $n_2=n+1$ for all $n\in\N$, we obtain
\Eq{*}{
\bigg(n\,,\,\dfrac{u_{n-1}+u_{n+1}}{2}\bigg)\in epi\big(u_n\big)\qquad \mbox{for all}\qquad n\in\N.
}
This provides $u_n\leq \dfrac{u_{n-1}+u_{n+1}}{2}$ and yields the convexity of the sequence $\big(u_n\big)_{n=0}^{\infty}$. This establishes the proposition.
\end{proof}

In the next result, we discuss the decomposition of any sequence through convex sequences.
\Prp{11}{
Every sequence can be represented as the difference of two convex sequences.
}
\begin {proof} 
Let $\big(u_n\big)_{n=0}^{\infty}$ be an arbitrary sequence. If $\big(u_n\big)_{n=0}^{\infty}$ is convex, there is nothing to prove. Hence, we assume that $\big(u_n\big)_{n=0}^{\infty}$ does not possess convexity. From a well-known result, this implies the existence of two monotonically increasing sequences 
$\big(v_n\big)_{n=0}^{\infty}$ and $\big(w_n\big)_{n=0}^{\infty}$ such that the expression $u_n:=v_n-w_n$ holds for all $n\in\N\cup\{0\}.$\\

Using these, we define two sequences $\big(x_n\big)_{n=0}^{\infty}$ and $\big(y_n\big)_{n=0}^{\infty}$ as follows:
\Eq{*}{
x_0=v_0 \quad \mbox{and}\quad x_n=&\sum_{i=0}^{n}v_i+\sum_{i=0}^{n-1}w_{i}\quad\mbox{for all}\quad n\in\N\\
&\qquad\,\,\mbox{and}\\
y_0=w_0 \quad \mbox{and}\quad y_n=&\sum_{i=0}^{n-1}v_i+\sum_{i=0}^{n}w_{i}\quad\mbox{for all}\quad n\in\N.
}
From the above construction, the following decomposition is evident
\Eq{*}{
u_n=v_n-w_n=x_n-y_n \quad  \mbox{ for all}\quad n\in\N\cup\{0\}.
}
Therefore, we only need to show the convexities of the sequences $\big(x_n\big)_{n=0}^{\infty}$ and $\big(y_n\big)_{n=0}^{\infty}$.\\

For all $n\in\N$, using monotonicity of the sequences $\big(v_n\big)_{n=0}^{\infty}$ and $\big(w_n\big)_{n=0}^{\infty}$, we can compute the following two inequalities
\Eq{*}{
x_n-x_{n-1}=v_n+w_{n-1}&\leq v_{n+1}+w_n=x_{n+1}-x_{n}\\
&\mbox{and}\\
y_n-y_{n-1}=w_n+v_{n-1}&\leq w_{n+1}+v_n=y_{n+1}-y_{n}.
}
The inequalities above yield that the sequences $\big(x_n\big)_{n=0}^{\infty}$ and $\big(y_n\big)_{n=0}^{\infty}$ carry usual convexities. This validates the proposition and completes the proof.
\end{proof}

A non-negative sequence is often referred to as $0$-convex, while increasingness is labelled as $1$-convex. The terms $2$-convex or ordinary convexity can be used interchangeably. Following this pattern, a sequence
$\big(u_n\big)_{n=0}^{\infty}$ is called \textit{$n$-convex} if for all $i\in\N\cup\{0\}$, it satisfies the following 
\Eq{*}{
u_{n+i}-\binom n 1 u_{n+i-1}+\cdots+(-1)^n  u_{i}\geq 0 \qquad \mbox{or equivalently}\qquad \sum_{k=0}^{n}(-1)^k\binom n k u_{n+i-k}\geq 0.
}
Using the constructive proofs above, one can generalize the sequence decomposition result to $n$-convexity.
\section{On Minorant, Sandwich and Stability results}
In this section, our primary objective is to propose a Hyers-Ulam-type stability result for convex sequences. We decide to break the proof into three parts. First, we demonstrate that under some minimal assumption, for every sequence that possesses a lower bound, it is viable to formulate a non-trivial convex minorant. Next, using the construction of the convex minorant, we derive a sandwich-type result, and finally we present a Hyers-Ulam-type stability result as a consequence of the previous findings.\\

A sequence $\big(v_n\big)_{n=0}^{\infty}$ is said to be a \textit{minorant} of $\big(u_n\big)_{n=0}^{\infty}$, if it satisfies the inequality $v_n\leq u_n$ for all $n\in\N\cup\{0\}.$ The study of minorants is crucial in approximation theory. In particular, nonconstant minorants that closely approximate a given sequence while satisfying additional regularity properties such as monotonicity, convexity, or subadditivity, etc. are especially interesting. Their construction reveals a subtle balance between approximation accuracy and structural constraints. In the next result, we present a method for constructing a nontrivial convex minorant of an arbitrary sequence.
\Prp{41}{Let $\big(u_n\big)_{n=0}^{\infty}$ be a sequence  bounded from below satisfying the condition 
\newline
$\inf\{u_n\}<u_0$, then there exists a non-trivial (non-constant) convex sequence $\big(v_n\big)_{n=0}^{\infty}$ such that the inequality $\inf\{u_n\}\leq v_n\leq u_n$ holds for all $n\in\N\cup\{0\}.$
}
\begin{proof}
Since the sequence $\big(u_n\big)_{n=0}^{\infty}$ is bounded below, the constant sequence $\Big(\inf\{u_n\}\Big)_{n=0}^{\infty}$ is a convex minorant of $\big(u_n\big)_{n=0}^{\infty}$. Our aim is to construct a sharper, non-constant minorant that possesses convexity.\\

To establish the result, we need to construct a class of sequences $\Bigg\{\Big(v_n^{^{k}}\Big)_{n=0}^{\infty}\,\bigg|\,k\in\N\cup\{0\}\Bigg\}$. 
\newline
First for $k=0$, we define the sequence $\Big(v_n^{^{0}}\Big)_{n=0}^{\infty}$ as follows
\Eq{42}{
v_{0}^{^{0}}=u_0\qquad 
\mbox{and}
\qquad
v_{n}^{^{0}}=\min
\Bigg\{
u_n\,,\,\dfrac{u_{n-1}+u_{n+1}}{2}
\Bigg\}
\quad \quad \mbox{for all} \quad n\in\N.
}
Using the above construction, for each $k\in\N$, we formulate the sequences $\Big(v_n^{^{k}}\Big)_{n=0}^{\infty}$ in a successive way as below 
\Eq{43}{
v_{0}^{^{k}}=u_0\qquad 
\mbox{and}
\qquad
v_{n}^{^{k}}=\min
\Bigg\{
v_n^{^{k-1}}\,,\, \dfrac{v_{n-1}^{^{k-1}}+v_{n+1}^{^{k-1}}}{2}
\Bigg\}
\quad \quad \mbox{for all} \quad n\in\N.
}
Now, from the $n^{th}$ elements of the sequences in the class  $\Bigg\{\Big(v_n^{^{k}}\Big)_{n=0}^{\infty}\,\bigg|\,k\in\N\cup\{0\}\Bigg\}$, we can extract a corresponding sequence $\Big(v_n^{^{k}}\Big)_{k=0}^{\infty}$. From our above construction, it is clear that the sequence $\Big(v_n^{^{k}}\Big)_{k=0}^{\infty}$ is monotonically decreasing. Besides, our boundedness assumption on $\big(u_n\big)_{n=0}^{\infty}$ yields
\Eq{22222}{
\inf\,\{u_n\}\leq v_n^{^{k}}\qquad \mbox{for all} \quad k\in\N.
}
This ensures that the sequence $\Big(v_n^{^{k}}\Big)_{n=0}^{\infty}$ is convergent. We assume 
$\underset{k\to\infty}{\lim}v_n^{^{k}}=v_n$. From our initial assumption $\inf\,\{u_n\}\leq u_0$, sequence construction in \eq{42}, \eq{43} and from the inequality \eq{22222}, it is evident that $\inf\,\{u_n\}\leq v_{n}$ holds for all $n\in\N\cup\{0\}$. To complete our theorem, it will be sufficient to show that the sequence 
$\big(v_n\big)_{n=0}^{\infty}$ is convex.\\

Now, we choose arbitrarily any $n^{th}$ element of the sequence $\Big(v_n^{^{k}}\Big)_{k=0}^{\infty}$. This yields
\Eq{*}{
v_n\leq v_n^{^{k}}\leq \dfrac{v_{n-1}^{^{k-1}}+v_{n+1}^{^{k-1}}}{2}\qquad \mbox{for all} \quad k\in\N.
}
Upon taking $\underset{k\to\infty}{\lim}$ in the right-most part of the above inequality, we get $v_n\leq \dfrac{v_{n-1}+v_{n+1}}{2}$. This establishes the convexity of the sequence $\big(v_n\big)_{n=0}^{\infty}$. Also, our assumption $\inf u_n<u_0$ ensures that the resultant convex sequence $\big(v_n\big)_{n=0}^{\infty}$ is non-constant. This completes the proof.
\end{proof}
In the above proposition, the condition $\inf\{u_n\}<u_0$ can be omitted. However, In certain cases, the convex minorant obtained in this way may be constant. One such example is the class of periodic sequences where $\inf\{u_n\}=u_0$ holds.\\

In function theory, one of the highly celebrated 
sandwich-type results that involves a convex function was proposed by  Baron, Matkowski, and Nikodem in their paper \cite{Baron}. The main result of their paper is the following theorem.
\Thm{10191}{
Let $f,g:I(\subseteq\R)\to\R$ are two real valued functions such that for all $x,y\in I$ and $t\in[0,1]$, the following functional inequality is satisfied
\Eq{*}{
f(tx+(1-t)y)\leq t\,g(x)+(1-t)\,g(y)\,;
}
then there exists a convex function $h:I\to\R$ such that $f(x)\leq h(x)\leq g(x)$ holds for all $x\in I$.
}
Through the following proposition, we try to achieve a similar result in a discrete setting. But first, we need to recall the definition of dyadic fractions. A rational number which can be represented as
$
\dfrac{m}{2^n},
$
where \(m \in \mathbb{Z}\) and \(n \in \mathbb{N}\cup\{0\}\) is called
\textit{dyadic rational number}. However, for our results, we are only interested in the set of all non-negative proper dyadic rational numbers, which can be defined as follows
\[
\mathbb{D}
=
\left\{
\dfrac{m}{2^n}
\;\middle|\;
m \in \mathbb{\N}\cup\{0\} \quad\mbox{and}\quad n \in \mathbb{N}\quad\mbox{with}\quad m\leq 2^n
\right\}.
\]
  
\Prp{44}{Let $\big(w_n\big)_{n=0}^{\infty}$ be a minorant of  the sequence $\big(u_n\big)_{n=0}^{\infty}$. Additionally, for all $n\in\N$, the sequence $\big(w_n\big)_{n=0}^{\infty}$ satisfies the following discrete functional inequality  
\Eq{45}{
w_n\leq\sum_{i=0}^{n-1}p_i\,u_i+q\,u_n+\sum_{i=n+1}^{\infty}r_i\,u_i\qquad\mbox{with}\qquad \sum_{i=0}^{n-1} p_i+q+\sum_{i=n+1}^\infty r_i=1\,\,;
}
where $p_i$'s, q and $r_i$'s $\in\mathbb{D}$ such that if some of the $p_i$'s are non-zero, then some $r_i$'s also attain positive value and vice versa. Then there exists a convex sequence $\big(v_n\big)_{n=0}^{\infty}$ such that $w_n\leq v_n\leq u_n$ holds for all $n\in\N\cup\{0\}.$
}
\begin{proof}
This proposition is a direct consequence of the above result. We start by constructing the class of sequences $\Bigg\{\Big(v_n^{^{k}}\Big)_{n=0}^{\infty}\, \bigg|\,k\in\N\cup\{0\}\Bigg\}$ as in \eq{42} and \eq{43}. From our construction, it is easily verifiable that for any $n^{th}$ element ($n\in\N$) of $k^{th}$ sequence ($k\in\N$), i.e. $v_n^{^{k}}$ can be represented as the following convex combination
\Eq{*}{
v_n^{^{k}}=\sum_{i=0}^{n-1}p_i\,u_i+q\,u_n+\sum_{i=n+1}^{\infty}r_i\,u_i \qquad\mbox{with}\qquad \sum_{i=0}^{n-1} p_i+q+\sum_{i=n+1}^\infty r_i=1;
}
where $p_i$'s, q and $r_i$'s $\in\mathbb{D}$ such that either simultaneously some of $p_i$'s and some of the $r_i$'s remain non-zero or instead all of $p_i$'s and $r_i$'s turn to be $0$. In view of \eq{45} and the assumed minorant condition, this yields that the sequences $\big(w_n\big)_{n=0}^{\infty}$ and $\Big(v_n^{^{k}}\Big)_{k=0}^{\infty}$ satisfy the following 
\Eq{*}{
w_n\leq v_n^{^{k}}\leq u_n \quad \mbox{for all}\quad k\in\N\cup\{0\} \qquad\qquad \Big(n\in\N\cup\{0\}\Big).
}
Implementing the same arguments and  methodology as in \prp{41}, we get the following
\Eq{*}{
w_n\leq \lim_{k\to\infty}v_n^{^{k}}\leq u_n \qquad \mbox{which implies}\quad w_n\leq v_n\leq u_n \qquad\mbox{for all}\quad n\in\N\cup\{0\}.
}
This validates the assertion and completes the proof of it.
\end{proof}
There are scopes of further investigation regarding the possibility of obtaining a similar sandwich-type result
by weakening the inequality \eq{45} of \prp{44}. However, this is a challenging task. For example, the most obvious form of such an inequality that seems worth studying instead of \eq{45} is the following 
\Eq{*}{
w_n\leq \dfrac{u_{n-1}+u_{n+1}}{2} \qquad \mbox{for all} \quad n\in\N.
}
Unfortunately, it turns out that this simple inequality alone is not sufficient to ensure the existence of a convex sequence between  $\big(u_n\big)_{n=0}^{\infty}$ and $\big(w_n\big)_{n=0}^{\infty}$. The following example is self-verifying and clarifies the argument. 
\Eq{*}{
u_n:=
\begin{cases}
n+1\qquad \mbox{if}\qquad n\leq 5\\
\quad 5 \qquad \mbox{for all}\quad n>5
\end{cases}
\qquad \mbox{and}\qquad
w_n:=
\begin{cases}
n\qquad \quad \mbox{if}\qquad n\leq 5\\
5 \qquad \mbox{for all}\quad n>5
\end{cases}.
}
\\

We now define the notion of approximately convex or 
$\varepsilon$-convex sequence. Let $\varepsilon>0$ be fixed. A sequence $\big(u_n\big)_{n=0}^{\infty}$ is said to be 
\textit{$\varepsilon$-convex} if for all $n\in\N$, it satisfies the following inequality 
\Eq{46}{
u_n\leq\sum_{i=0}^{n-1}p_i\,u_i+q\,u_n+\sum_{i=n+1}^{\infty}r_i\,u_i+\varepsilon\qquad\mbox{with}\qquad \sum_{i=0}^{n-1} p_i+q+\sum_{i=n+1}^\infty r_i=1\,\,;
}
where $p_i$'s, q and $r_i$'s $\in \mathbb{D}$ such that if some of $p_i$'s are non-zero then some $r_i$'s also positive and vice versa. \\

We show that an $\varepsilon$-convex sequence can be closely approximated by a convex sequence. Also, for our next result, the terminology of supremum norm as the difference of two real-valued sequences, $\big(u_n\big)_{n=0}^{\infty}$ and $\big(v_n\big)_{n=0}^{\infty}$ is required, which can be defined as follows
\Eq{*}{
\|v_n-u_n\|_{\infty}=\sup_{n\in\N\cup\{0\}}\bigg\{|u_n-v_n|\bigg\}.
}
Now, we can present the main result of this section.
\Thm{47}{If $\big(u_n\big)_{n=0}^{\infty}$ is an 
$\varepsilon$-convex sequence, then there exists a convex sequence 
$\big(v_n\big)_{n=0}^{\infty}$ such that the inequality $\|u_n-v_n\|_{\infty}\leq \varepsilon/2$ holds. Conversely, if a convex sequence $\big(v_n\big)_{n=0}^{\infty}$ satisfies the inequality $\|u_n-v_n\|_{\infty}\leq \varepsilon/2$ , then the sequence $\big(u_n\big)_{n=0}^{\infty}$ is 
$\varepsilon$-convex.}
\begin{proof}
Since the sequence $\big(u_n\big)_{n=0}^{\infty}$ is approximately convex, the inequality \eq{46} holds along with the mentioned coefficient-related conditions.
We can rewrite this inequality as follows
\Eq{*}{
u_n-\varepsilon\leq\sum_{i=0}^{n-1}p_i\,u_i+q\,u_n+\sum_{i=n+1}^{\infty}r_i\,u_i\qquad (n\in\N).
}
Now, considering $w_n=u_n-\varepsilon$ for all 
$n\in\N\cup\{0\}$, we have the sequence $\big(w_n\big)_{n=0}^{\infty}$ 
that satisfies the inequality \eq{45} of \prp{44} along with the mentioned conditions on $p_i$'s, $q$ and $r_i$'s. In view of the \prp{44}, we have a convex sequence 
$\big(\widetilde{v}_n\big)_{n=0}^{\infty}$ such that the following inequality holds
\Eq{48}{
u_n-{\varepsilon}\leq \widetilde{v}_n\leq u_n \qquad \mbox{for all} \quad n\in\N\cup\{0\}.
}
Shifting the sequence $\big(\widetilde{v}_n\big)_{n=0}^{\infty}$ by $\varepsilon/2$, we still obtain a convex sequence $\big({v}_n\big)_{n=0}^{\infty}$. Thus using this translation, the inequality \eq{48}, can be re-structured as below
\Eq{49}{
u_n-\dfrac{\varepsilon}{2}\leq {v}_n\leq u_n+\dfrac{\varepsilon}{2} \qquad \mbox{for all} \quad n\in\N\cup\{0\}.
}
This validates the first assertion.\\

To show the converse part, we assume that sequence $\big({v}_n\big)_{n=0}^{\infty}$ possesses ordinary convexity and \eq{49} holds. Then from the left-most inequality of \eq{49}, we can compute the following 
\Eq{*}{
u_n\leq v_n+\dfrac{\varepsilon}{2}\leq \dfrac{v_{n-1}+v_{n+1}}{2}+\dfrac{\varepsilon}{2}\qquad (n\in\N)\,.
} 
Now we can keep extending the right-most part of the above inequalities simply by using the convexity of $\big({v}_n\big)_{n=0}^{\infty}$. At any instance, this extension appears to be in the following form 
\Eq{*}{
u_n\leq\Bigg(\sum_{i=0}^{n-1}p_i\,v_i+q\,v_n+\sum_{i=n+1}^{\infty}r_i\,v_i\Bigg)+\dfrac{\varepsilon}{2}\qquad\mbox{with}\qquad \sum_{i=0}^{n-1} p_i+q+\sum_{i=n+1}^\infty r_i=1\,\,;
}
where $p_i$'s, q and $r_i$'s $\in\mathbb{D}$ such that simultaneously some $p_i$'s and some $r_i$'s remain positive. Using the right-most inequality of \eq{49}, the above inequality can be extended as follows
\begin{small}
\Eq{*}{
u_n\leq\Bigg(\sum_{i=0}^{n-1}p_i\,\bigg(u_i+\dfrac{\varepsilon}{2}\bigg)+q\,\bigg(u_n+\dfrac{\varepsilon}{2}\bigg)+\sum_{i=n+1}^{\infty}r_i\,\bigg(u_i+\dfrac{\varepsilon}{2}\bigg)\Bigg)+\dfrac{\varepsilon}{2}\quad\mbox{with}\quad \sum_{i=0}^{n-1} p_i+q+\sum_{i=n+1}^\infty r_i=1\,.
}  
\end{small}
Upon simplification, the above inequality leads to \eq{46} which implies that the sequence $\big({u}_n\big)_{n=0}^{\infty}$ possesses $\varepsilon$-convexity. This completes the proof of the statement.
\end{proof}
The rest of the paper is devoted to the analysis of subsequences.
\section{On Convex Subsequences}
The motivation behind this section has already been discussed in detail in the introduction. 
To present our main result in this section, we first require an equivalent definition of a convex function. The slope-oriented definition of an ordinary convex function was introduced in Theorem 3.7.1 of the book \cite{Kuczma}. For our use, the result can be stated in a simplified form as follows
\Thm{599}{
A function $f:\R_+\to\R$ is convex if and only if for all $x,y,z\in\R_+$ with $x<y<z$, the following functional inequality holds
\Eq{122}{
\dfrac{f(y)-f(x)}{y-x}\leq \dfrac{f(z)-f(y)}{z-y}.
}
}
We also recall the definition of the sign function, which extracts the sign of a real number. In our context, the sign function, $sgn:\N\cup\{0\}\to\{-1,0,1\}$ is defined as follows
\Eq{*}{
sgn(\ell-k)=
\begin{cases}
\,\,\,\,1\,\,\,\,\,\,\,\,\,\quad\mbox{if}\quad k<\ell\\
\,\,\,\,0\,\,\,\,\,\,\,\,\,\quad\mbox{if}\quad k=\ell\\
-1\,\,\,\mbox{else if}\quad k>\ell
\end{cases}.
}
\Thm{601}{Let $\big(u_n\big)_{n=0}^{\infty}$ be a concave sequence and $\Big(u_{n_{_{i}}}\Big)_{i=0}^{\infty}$ is a subsequence of it. Then the following assertions are equivalent to each other.
\begin{enumerate}[(i)]
\item  $\Big(u_{n_{_{i}}}\Big)_{i=0}^{\infty}$ is convex.
\vspace{3.5mm}
\item For all $k,\ell,m\in\N\cup\{0\}$ with $k<\ell<m$, the following inequality holds
\Eq{602}{
\bigg(\dfrac{1}{\ell-k}-\dfrac{1}{n_\ell-n_k}\bigg)\Big(u_{n_{_{\ell}}}-u_{n_{_{k}}}\Big)
\leq 
\bigg(\dfrac{1}{m-\ell}-\dfrac{1}{n_m-n_\ell}\bigg)
\Big(u_{n_{_{m}}}-u_{n_{_{\ell}}}\Big).
}
\vspace{2mm}
\item There exists a monotone sequence $\big(v_{n_{_{\ell}}}\big)_{\ell=0}^{\infty}$ such that for all $k,\ell\in\N\cup\{0\}$,  the following discrete functional inequality is satisfied
\Eq{603}{
\Big(\left|n_k-n_{\ell}\right|-|k-\ell|\Big)\big(u_{n_{_{\ell}}}-u_{n_{_{k}}}\big)\leq 
sgn(\ell-k)\,v_{n_{_{\ell}}} (\ell-k)(n_{\ell}-n_k) \,.
}
\end{enumerate}
}
\begin{proof}
$(i)\Rightarrow (ii) :\,$ Let $\big(u_{n_{_{i}}}\big)_{i=0}^{\infty}$ be a convex subsequence of the concave sequence   $\big(u_n\big)_{n=0}^{\infty}$. First, using the \thm{15} and then applying inequality  \eq{122} of \thm{599}, for all $k,\ell,m\in\N\cup\{0\}$ with $k<\ell<m$, we can conclude the following two inequalities
\Eq{*}{
\dfrac{u_{n_{_{\ell}}}-u_{n_{_{k}}}}{\ell-k}
\leq 
\dfrac{u_{n_{_{m}}}-u_{n_{_{\ell}}}}{m-\ell}
\qquad \mbox{and}\qquad 
-\dfrac{u_{n_{_{\ell}}}-u_{n_{_{k}}}}{n_\ell-n_k}
\leq 
-\dfrac{u_{n_{_{m}}}-u_{n_{_{\ell}}}}{n_m-n_\ell}.
}
Summing up the two inequalities above side by side, we arrive at the inequality \eq{602}. This establishes the assertion $(ii).$\\

$(ii)\Rightarrow (iii):\,$ We assume that the assertion $(ii)$ holds. Now, we consider $\ell\in\N\cup\{0\}$ be fixed and define the sequence $\Big(v_{n_{_{\ell}}}\Big)_{\ell=0}^{\infty}$ as follows:
\Eq{202020}{
v_{n_{_{\ell}}}:=\inf_{\ell<m}\,\bigg(\dfrac{1}{m-\ell}-\dfrac{1}{n_{m}-n_{\ell}}\bigg)\Big(u_{n_{_{m}}}-u_{n_{_{\ell}}}\Big)\qquad \quad \Big(m\in\N\Big).
}
Using this, the inequality \eq{602} can be refined as follows
\Eq{605}{
\bigg(\dfrac{1}{\ell-k}-\dfrac{1}{n_\ell-n_k}\bigg)\Big(u_{n_{_{\ell}}}-u_{n_{_{k}}}\Big)
\leq v_{n_{_{\ell}}}\leq
\bigg(\dfrac{1}{m-\ell}-\dfrac{1}{n_m-n_\ell}\bigg)
\Big(u_{n_{_{m}}}-u_{n_{_{\ell}}}\Big) \quad \big(k<\ell<m\big).
}
From the left-most inequality of \eq{605}, we obtain the following inequality
\Eq{606}{
\Big(\big(n_{\ell}-n_k\big)-(\ell-k)\Big)\big(u_{n_{_{\ell}}}-u_{n_{_{k}}}\big)\leq 
 v_{n_{_{\ell}}} (\ell-k)(n_{\ell}-n_k) \qquad \qquad
 \big(k,\,\ell\in\N\cup\{0\}\quad\mbox{with}\quad k<{\ell}\big).
}
Similarly, from the right-most inequality of \eq{605}, we get the inequality below
\Eq{*}{
\Big(\big(n_m-n_{\ell}\big)-(m-\ell)\Big)\big(u_{n_{_{\ell}}}-u_{n_{_{m}}}\big)\leq 
-v_{n_{_{\ell}}} (\ell-m)(n_{\ell}-n_m) \qquad \qquad
\big(\ell,\,m\in\N\quad\mbox{with}\quad {\ell}<m\big).
}
From the definition of $\Big(v_{n_{_{\ell}}}\Big)_{\ell=0}^{\infty}$ in \eq{202020}, the inequality above can be extended to all $\ell,\,m\in\N\cup\{0\}$, and we get the following
\Eq{6066}{
\Big(\big(n_m-n_{\ell}\big)-(m-\ell)\Big)\big(u_{n_{_{\ell}}}-u_{n_{_{m}}}\big)\leq 
-v_{n_{_{\ell}}} (\ell-m)(n_{\ell}-n_m) \quad 
\big(\ell,\,m\in\N\cup\{0\}\quad\mbox{with}\quad {\ell}<m\big).
}
Replacing the variable $m$ with $k$ in the above inequality, we have
\Eq{607}{
\Big(\big(n_k-n_{\ell}\big)-(k-\ell)\Big)\big(u_{n_{_{\ell}}}-u_{n_{_{k}}}\big)\leq 
-v_{n_{_{\ell}}} (\ell-k)(n_{\ell}-n_k) \qquad 
\big(k,\,\ell\in\N\cup\{0\}\quad\mbox{with}\quad {\ell}<k\big).
}
The inequalities \eq{606} and \eq{607}, together yields the following
\begin{small}
\Eq{*}{
\Big(\left|n_{\ell}-n_k\right|-|\ell-k|\Big)\big(u_{n_{_{\ell}}}-u_{n_{_{k}}}\big)\leq sgn(\ell-k)\, 
v_{n_{_{\ell}}} (\ell-k)(n_{\ell}-n_k) \qquad 
\big(k,\,\ell\in\N\cup\{0\}\quad\mbox{with} \quad k\neq {\ell}\big). 
}
\end{small}
Also, one can easily observe that the above inequality holds as equality for the case $k= {\ell}$. This establishes the inequality \eq{603}. Now, we are only required to show the monotonicity(increasingness) of the sequence $\Big(v_{n_{_{\ell}}}\Big)_{\ell=0}^{\infty}$.\\

Let $k,\ell\in\N\cup\{0\}$ with $k<\ell$. Clearly, this also indicates $n_{_{k}}<n_{_{\ell}}$. From the validated inequality \eq{603}, we have \eq{606} and the following inequality 
\Eq{610}{
\Big(\big(n_{\ell}-n_k\big)-(\ell-k)\Big)\big(u_{n_{_{k}}}-u_{n_{_{\ell}}}\big)\leq 
-v_{n_{_{k}}} (l-k)(n_{\ell}-n_k) \qquad \quad
 \big(k<{\ell}\big).
}
Adding up the inequalities \eq{606} and \eq{610} side by side, we arrive at the following
\Eq{*}{
0\leq \Big(v_{n_{_{\ell}}}-v_{n_{_{k}}}\Big)(\ell-k)(n_{\ell}-n_k).
}
Our assumption regarding $k,\ell$ gives us $(\ell-k)(n_{\ell}-n_k)>0$. Thus, from the above inequality, we conclude that $v_{n_{_{k}}}\leq v_{n_{_{\ell}}}$ holds. Since $k,\ell$ were arbitrarily chosen, this implies that the sequence 
$\Big(v_{n_{_{\ell}}}\Big)_{\ell=0}^{\infty}$ is increasing. This proves the assertion $(iii)$.\\

$(iii)\Rightarrow (i):\,$ We assume that the assertion $(iii)$ holds. Let $k,\ell,m\in\N\cup\{0\}$ such that $k<\ell<m$. This also ensures the inequality $n_k<n_{\ell}<n_m$. Using, the inequality \eq{603} of the assertion $(iii)$, we can obtain the inequalities \eq{606} and \eq{6066}. First, restructuring and then combining these two inequalities, we arrive at the following
\Eq{909}{
\Bigg(\dfrac{1}{\ell-k}-\dfrac{1}{n_{\ell}-n_k}\Bigg)\big(u_{n_{_{\ell}}}-u_{n_{_{k}}}\big)\leq 
 v_{n_{_{\ell}}}\leq \Bigg(\dfrac{1}{m-\ell}-\dfrac{1}{n_{m}-n_\ell}\Bigg)\big(u_{n_{_{m}}}-u_{n_{_{\ell}}}\big).
}
By using the definition of $v_{n_{_{\ell}}}$, from the left-most inequality of \eq{909}, we can conclude the following
\Eq{10000}{
\dfrac{u_{n_{_{\ell}}}-u_{n_{_{k}}}}{\ell-k}\leq  v_{n_{_{\ell}}}+\dfrac{u_{n_{_{\ell}}}-u_{n_{_{k}}}}{n_\ell-n_k}\leq \dfrac{u_{n_{_{m}}}-u_{n_{_{\ell}}}}{n_m-n_\ell}+
\dfrac{u_{n_{_{\ell}}}-u_{n_{_{k}}}}{n_\ell-n_k}.
}
From the assertion $(ii)$, we can also define the non-decreasing sequence $\Big(w_{n_{_{\ell}}}\Big)_{\ell=1}^{\infty}$ as follows:
\Eq{*}{
w_{n_{_{\ell}}}:=\sup_{k<\ell}\,\bigg(\dfrac{1}{\ell-k}-\dfrac{1}{n_\ell-n_k}\bigg)\Big(u_{n_{_{\ell}}}-u_{n_{_{k}}}\Big)\qquad \quad \Big(k\in\N\cup\{0\}\Big).
}
Following the same methodology and similar mathematical steps, we can conclude the following inequality
\Eq{910}{
\Bigg(\dfrac{1}{\ell-k}-\dfrac{1}{n_{\ell}-n_k}\Bigg)\big(u_{n_{_{\ell}}}-u_{n_{_{k}}}\big)\leq 
w_{n_{_{\ell}}}\leq \Bigg(\dfrac{1}{m-\ell}-\dfrac{1}{n_{m}-n_\ell}\Bigg)\big(u_{n_{_{m}}}-u_{n_{_{\ell}}}\big).
}
From the right-most inequality of \eq{910}, we can obtain the following extended inequality
\Eq{10001}{ 
 \dfrac{u_{n_{_{m}}}-u_{n_{_{\ell}}}}{n_m-n_\ell}+
\dfrac{u_{n_{_{\ell}}}-u_{n_{_{k}}}}{n_\ell-n_k}
\leq w_{n_{_{\ell}}}+\dfrac{u_{n_{_{m}}}-u_{n_{_{\ell}}}}{n_{m}-n_\ell}\leq\dfrac{u_{n_{_{m}}}-u_{n_{_{\ell}}}}{m-\ell}.
}
Combining the inequalities, \eq{10000} and \eq{10001}, we observe the following 
\Eq{*}{
\dfrac{u_{n_{_{\ell}}}-u_{n_{_{k}}}}{\ell-k}\leq\dfrac{u_{n_{_{m}}}-u_{n_{_{\ell}}}}{m-\ell}
\quad \mbox{for all}\quad k,\ell,m\in\N\cup\{0\}\quad \mbox{with}\quad k<\ell<m.
}
This yields that the subsequence $\big(u_{n_{_{i}}}\big)_{i=0}^{\infty}$ is convex, i.e. assertion $(i)$. This completes our characterization. 
\end{proof} 
To present our next proposition, we need to define a new concept, called a convex triplet. Any triplet $n_{_{k-1}}$, $n_{_{k}}$, and $n_{_{k+1}}\in\N\cup\{0\}$ with $n_{_{k-1}}<n_{_{k}}<n_{_{k+1}}$ is called a \textit{convex triplet}, if it satisfies the following inequality condition 
\Eq{801}{
0\leq n_{_{k}}-n_{_{k-1}}\leq n_{_{k+1}}-n_{_{k}}.
} 
In the next proposition, we show that if a sequence is increasing and convex, then any subsequence whose domain is sequentially convex in $\N\cup\{0\}$ also possesses convexity. 
\Prp{800}{Let $\big(u_n\big)_{n=0}^{\infty}$ be a real valued   non-decreasing convex sequence. If for all $k\in\N$, the triplet $n_{_{k}}, n_{_{k-1}}$ and $n_{_{k+1}}\in\N$ is convex, then the subsequence $\Big(u_{n_{k}}\Big)_{k=0}^{\infty}$ possesses convexity.}
\begin{proof}
Let $k\in\N$ be arbitrary and $n_{_{k-1}}$, $n_{_{k}}$, and $n_{_{k+1}}\in\N$ form a convex triplet. Therefore, it satisfies the inequality \eq{801}.
Also, as discussed in the previous result, the convexity along with increasingness of the sequence $\big(u_n\big)_{n=0}^{\infty}$ provides the following
\Eq{802}{
0\leq \dfrac{u_{n_{_{k}}}-u_{n_{_{k-1}}}}{n_{_{k}}-n_{_{k-1}}}\leq \dfrac{u_{n_{_{k+1}}}-u_{n_{_{k}}}}{n_{_{k+1}}-n_{_{k}}}.
}
Multiplying the right-most inequalities of \eq{801} and \eq{802} side by side, we arrive at
\Eq{*}{
u_{n_{_{k}}}-u_{n_{_{k-1}}}\leq u_{n_{_{k+1}}}-u_{n_{_{k}}}.
}
This demonstrates the subsequence $\Big(u_{n_{k}}\Big)_{k=0}^{\infty}$ is convex and completes the proof.
\end{proof}
A subsequence $\Big(u_{n_{k}}\Big)_{k=0}^{\infty}$ of the sequence $\big(u_n\big)_{n=0}^{\infty}$ is said to be \textit{periodically indexed subsequence}, if there exists a positive constant $\L$ such that for each $k\in\N$, the equality $n_k-n_{k-1}=\L$ holds. We can also represent such a sequence with the notation
$\big(u_{{_{i+k\L}}}\big)_{k=0}^{\infty}$, where $u_{_{i}}$ is the initial term, $\L$ is the common index difference. \\

A real valued sequence $\big(u_n\big)_{n=0}^{\infty}$ is called \textit{subadditive} if  the following discrete inequality holds
\Eq{*}{
u_{n_1+n_2}\leq u_{n_1}+u_{n_2}\qquad \mbox{for all}\qquad n_1,n_2\in\N. 
}
Many classes of sequences can be studied within the broader framework of subadditive sequences. Several results demonstrate that, under certain conditions, monotonicity, convexity, and subadditivity are closely related in the discrete setting. For instance, every positive decreasing sequence is subadditive. Moreover, every non-negative concave sequence is also subadditive. In the following result, we show how convexity can be used to refine the subadditivity condition.
\Prp{810}{
Let $\big(u_n\big)_{n=0}^{\infty}$ be a subadditive sequence and $\big(u_{{_{i+k\L}}}\big)_{k=0}^{\infty}$ is a periodically indexed subsequence of it. If $\big(u_{_{{i+k\L}}}\big)_{k=0}^{\infty}$ is convex, then for all $k\in\N$, the following functional inequality holds
\Eq{811}{
u_{_{2({i+k\L})}}\leq 2u_{{_{i+k\L}}
}\leq u_{{_{i+(k-1)\L}}}+u_{{_{i+(k+1)\L}}}.
}
}
\begin{proof}
Since $\big(u_n\big)_{n=0}^{\infty}$ possesses subadditivity, the following two inequalities are obvious
\Eq{812}{
u_{_{2(i+k\L)}}&=u_{_{{(i+(k-1)\L)}+{(i+(k+1)\L)}}} \leq u_{{_{i+(k-1)\L}}}+u_{{_{i+(k+1)\L}}}\\
&\qquad\qquad\qquad\mbox{and}\\
&\quad u_{_{2(i+k\L)}}=u_{{_{i+k\L}}+{{i+k\L}}} \leq 2u_{{_{i+k\L}}
}.
}
Therefore, the inequality \eq{811} provides ordering to all terms that form the two inequalities above.  
The convexity of the subsequence $\big(u_{n_{{i+k\L}}}\big)_{k=0}^{\infty}$ yields the following inequality
\Eq{*}{
2u_{{_{i+k\L}}
}\leq u_{{_{i+(k-1)\L}}}+u_{{_{i+(k+1)\L}}}.
} 
This inequality together with the second inequality of \eq{812} provides \eq{811}. This completes the proof.
\end{proof}
The research presented in this paper offers numerous opportunities for further exploration. One obvious direction is to investigate stability results and subsequence characterizations for higher-order discrete convexity.  Moreover, similar to \thm{47}, one may consider different variants of generalized convex sequences and establish corresponding stability results. Analogous to \thm{601}, it would also be interesting to characterize subadditive sequences that contain a superadditive subsequence.

\section*{Acknowledgements}

The author sincerely thanks \textbf{István Szalkai} for his valuable discussions, constructive suggestions, and critical remarks, which significantly improved the quality and clarity of this manuscript.
\section*{Statements and Declarations}
\noindent\textbf{Funding.} 

\vspace{6pt}
\noindent\textbf{Competing Interests.} 
The author declare that there are no financial or non-financial competing interests relevant to the contents of this article.

\vspace{6pt}
\noindent\textbf{Ethics Approval.} 
Not applicable. This study does not involve human participants or animals.

\vspace{6pt}
\noindent\textbf{Consent to Participate.} 
Not applicable.

\vspace{6pt}
\noindent\textbf{Consent for Publication.} 
Not applicable.

\vspace{6pt}
\noindent\textbf{Data, Materials and/or Code Availability.} 
No datasets or code were generated or analysed during the current study.

\end{document}